\documentclass[a4paper,12pt]{article}
\usepackage{amsfonts, amsbsy}
\usepackage{graphics}
\usepackage{amsmath}
\usepackage{latexsym}
\usepackage{amssymb}

\newcommand{\spa}{\vspace{.2in}}
\newcommand{\bo}{\textbf}
\newcommand{\noi}{\noindent}
\newcommand{\LA}{\mathcal{L}}
\newcommand{\U}{\mathcal{U}}
\newcommand{\Rd}{\mathbb{R}^d}
\newcommand{\eps}{\epsilon}
\newcommand{\defn}{\stackrel{def}{=}}
\date{}
\begin{document}
\begin{center}\bo{\large SMALL NOISE ASYMPTOTICS FOR INVARIANT\\ \ \\
DENSITIES FOR A CLASS OF DIFFUSIONS:\\ \ \\
A CONTROL THEORETIC VIEW}
\end{center}

\vspace{.5in}

\begin{center}
ANUP BISWAS\footnote{Centre for Applicable
Mathematics, Tata Institute of Fundamental Research,
Post Bag No. 03, Sharadanagar, Chikkabommasandra, Bangalore 560065, India}
AND VIVEK S. BORKAR\footnote{School of Technology and Computer Science,
Tata Institute of Fundamental Research, Homi Bhabha Rd., Mumbai 400005,
India. (email: borkar@tifr.res.in) Work supported in part by a J. C.
Bose Fellowship.}
\end{center}

\vspace{.5in}

\noindent \bo{Abstract:}
 We consider multidimensional nondegenerate diffusions with invariant densities,
 with the diffusion matrix scaled by a small $\eps > 0$. The o.d.e. limit
corresponding to $\eps = 0$ is assumed to have the origin as its unique globally
asymptotically stable equilibrium. Using control theoretic methods, we show
that in the $\eps \downarrow 0$ limit, the invariant density has the form
$\approx exp(-W (x)/\eps^ 2 )$,
where the $W$ is characterized as the optimal cost of a deterministic control
problem. This generalizes an earlier work of Sheu. Extension to multiple
equilibria is also given.

\spa


\section{Introduction}
A recurrent theme in applied mathematics is the resolution of non-uniqueness
issue of a deterministic system by considering its perturbation with small
noise and recovering a hopefully unique choice by passing to the vanishing
noise limit. (This idea is attributed to Kolmogorov in \cite{EcRu}, p. 626)
To
mention two such instances, see the analysis of jump phenomena in nonlinear
circuits in \cite{Sas} and equilibrium selection in evolutionary games in \cite{FoYo}.
 In fact
the notion of `viscosity solutions' we use later in this work can be motivated
along these lines.

\spa

    A similar situation also arises in the passage from quantum to classical
mechanics (the correspondence principle) when one uses the Feynman-Kac
device to convert Feynman path integrals to Wiener space integrals by complexifying time, whence the problem reduces to that of small noise behavior
of diffusions. Motivated by this, Freidlin and Wentzell have extensively de-
veloped the small noise limit theory for diffusions \cite{Fre}.

\spa

    Building on the classic work of Freidlin and Wentzell \cite{Fre}, Sheu in \cite{She}
characterized the small noise limit of the invariant density of a positive recur
rent diffusion for a restricted class of diffusions. This was further extended
by Day \cite{Day1}. Our aim is to establish a similar result in a more general set-up,
using a novel control theoretic approach based on the theory of viscosity solutions. We do so by exploiting the fact that the `adjoint' partial differential
equation satisfied by the invariant density gets converted to the Hamilton-Jacobi-Bellman equation of an associated ergodic control problem via the
logarithmic transformation. The limiting case thereof as the noise decreases
to zero can then be handled by exploiting the machinery of viscosity solutions.

\spa

    The paper is organized as follows: The next section introduces the problem and states the main result. Section $3$ uses the equivalent control formulation in order to go to the small noise limit and obtain the asymptotic
expression for the invariant density as $\approx exp(-W (x)/\eps^2 )$. Section
$4$ obtains
a representation for $W$ in terms of a deterministic control problem. Section
$5$ extends the results to the case of multiple equilibria.

\spa

    See also \cite{Ari1}, \cite{Ari2}, \cite{Mik}, \cite{stet} for related work.

\section{The small noise asymptotics}
We consider an $d$-dimensional diffusion $X$ given by the stochastic differential
equation
\begin{equation}
dX(t) = b(X(t))dt + \eps\sigma(X(t))dB(t), \ t \geq 0, \label{sde}
\end{equation}
where $b: \Rd \rightarrow \Rd, \sigma: \Rd \rightarrow \mathbb{R}^{d\times m}$ are smooth,
bounded, with bounded derivatives, and $B$ is an
$m$-dimensional Brownian motion. We assume:
\begin{enumerate}
\item $\|D^{\kappa}b(x)\|, \|D^{\kappa}\sigma(x)\| \leq C < \infty$ for some $C$ and all $x$
and all multi-indices $\kappa$.

\item $a(x)= [[a_{ij}(x)]] \defn \sigma(x)\sigma^T(x)$ satisfies:
\begin{displaymath}
\lambda\|x\|^2 \leq x^Ta(y)x \leq \Lambda\|x\|^2, \ \forall x, y\in\Rd,
\end{displaymath}
for some $0 < \lambda \leq \Lambda$.

\item There exist constants $\alpha > 0, 0 < \beta \leq 1$ such that
\begin{equation}
\limsup_{\|x\|\rightarrow\infty}\Big[\alpha\Lambda + \|x\|^{1 -
\beta}b(x)^T\left(\frac{x}{\|x\|}\right)\Big] < 0. \label{assum}
\end{equation}

\item $b(0) = 0$ and zero is the unique globally asymptotically stable equilibrium point
of the o.d.e
\begin{equation}
\dot{x}(t) = b(x(t)). \label{ode}
\end{equation}
\end{enumerate}

Let
\begin{displaymath}
\LA^{\eps} \defn \langle b(\cdot), \nabla \rangle + \frac{\eps^2}{2}\mbox{tr}\left(a(\cdot)\nabla^2\right)
\end{displaymath}

\spa

\noi {\bf Theorem 1 :}
 For each $\eps \in (0, 1)$, $X(\cdot)$ is
positive recurrent and the corresponding (unique) stationary
distributions $\{\mu^{\eps}, 0 < \eps < 1\}$ are tight. Furthermore,
$\lim_{\eps\downarrow 0}\mu^{\eps} = \delta_0$, where $\delta_x$ for
$x \in \Rd$ is the Dirac measure at $x$.

\spa

\noi \bo{Proof:} For $\mathcal{V}(x) \defn \|x\|^2$, direct computation shows that
\begin{eqnarray}
\LA^{\eps}\mathcal{V}(x) &=& \eps^2\sum_ia_{ii}(x) + 2b(x)^Tx \nonumber \\
&\leq& d\eps^2\Lambda + 2b(x)^Tx\ \ \nonumber \\
&\leq& d\eps^2\Lambda - \alpha\Lambda\|x\|^{\beta} \label{th1dis}
\end{eqnarray}
for $\|x\| \geq$ a sufficiently large $K > 0$. (The last inequality
follows by virtue of (\ref{assum}).) In particular,
\begin{displaymath}
\LA^{\eps}\mathcal{V}(x) < - \delta
\end{displaymath}
for some $\delta > 0$ and sufficiently large $\|x\|$. But this is
the standard stochastic Liapunov condition for positive recurrence
of nondegenerate diffusions, implying in particular the existence of
$\mu^{\eps}$. To see this, first note that a simple application of Dynkin formula in conjunction to the above leads to the conclusion that the mean hitting time of $\{x : \|x\| \leq K\}$ from any point is finite (p.\ 305, \cite{GihSko}). An invariant probability measure can then be constructed as in \cite{Khas} (see also \cite{Khas2}). Uniqueness follows from the observation that under the
stated hypotheses, $X(t)$ has a density $p( y | x, t) > 0$ for all
$t$ and
therefore so does $\mu^{\eps}(dx) = \int p(x|t, z)\mu^{\eps}(dz)$.
Thus any two invariant measures would have to be mutually absolutely
continuous. But by a well known fact from ergodic theory of Markov
processes (section XIII.4, \cite{Yosi}), two distinct extremal invariant probability measures of a
Markov process must be mutually singular. This contradiction implies
that the invariant probability  measure is unique. Also, by
Proposition 2.3 of
\cite{Met}, we have
\begin{equation}
\int |\LA^{\eps}\mathcal{V}|d\mu^{\eps} \leq K < \infty \label{metabdd}
\end{equation}
for some $K$. By Theorem 9.17 of \cite{Eth}, $\mu^{\eps}$ is
characterized by
\begin{equation}
\int \LA^{\eps}fd\mu^{\eps} = 0 \label{Eche}
\end{equation}
for twice continuously differentiable $f$ with bounded first and
second derivatives.

\spa

From (\ref{th1dis}), $\LA^{\eps}\mathcal{V}(x) \leq -g(x)$ for a $g$ satisfying
$\displaystyle{\lim_{\|x\|\rightarrow\infty}}g(x) = \infty$. Thus for a suitable $K^* > 0$,
$\LA^{\eps}\mathcal{V}(x) < 0$ for $\|x\| \geq K^*$. Using (\ref{metabdd}), we then
have
\begin{displaymath}
\int_{\|x\| \geq K^*} g(x) d\mu^{\eps} \leq \int_{\|x\| \leq K^*}|\LA^{\eps}\mathcal{V}|d\mu^{\eps}
\leq K.
\end{displaymath}
 Since
$\displaystyle{\lim_{\|x\|\rightarrow\infty}}g(x) = \infty$, this
implies tightness of $\{\mu^{\eps}, \eps > 0\}$. Now let
$\eps\downarrow 0$ in (\ref{Eche}) to obtain
\begin{displaymath}
\int \LA^0 fd\mu = 0
\end{displaymath}
for all $f$ as above and all limit points $\mu$ in $\mathcal{P}(\Rd)$ of
$\mu^{\eps}$ as $\eps\downarrow 0$. By the criterion (\ref{Eche}),
$\mu$ must be an  invariant measure for (\ref{ode}). But under our
assumptions on (\ref{ode}), this is possible only for $\mu =
\delta_0$. \hfill $\Box$

\spa

Our aim is to capture the precise manner in which $\mu^{\eps} \rightarrow
\delta_0$. Note that for $\eps > 0$, $\mu^{\eps}$ will have a
density $\varphi^{\eps}$. First we will prove the following result

\spa

\noi \bo{Theorem 2 :} $\lim_{\eps\downarrow
0}\eps^2\ln(\varphi^{\eps}(x)) = - W(x)$ for
\begin{equation}
W(x) \defn \inf\Big(\int_0^{\infty}u(t)^Ta^{-1}(y(t))u(t)dt\Big),
\label{Repr}
\end{equation}
where the infimum is over all measurable $u(\cdot)$ such that the trajectory
$y(\cdot)$ of
\begin{eqnarray}
 \dot{y}(t)=-b(y(t)) - u(t), \ y(0)=x, \label{controls}
\end{eqnarray}
satisfies: $y(t)\stackrel{t\uparrow\infty}{\to}0$.

\spa

We prove this in the subsequent sections.

\section{The control formulation}

As observed above
\begin{equation}
\int\LA^{\eps}fd\mu^{\eps} = \int \LA^{\eps}f(x)\varphi^{\eps}(x)dx = 0
\label{inveq}
\end{equation}
for smooth compactly supported $f$. Hence $\varphi^{\eps}$ satisfies
\begin{equation}
\frac{\eps^2}{2}\sum_{i,j=1}^na_{ij}\varphi^{\eps}_{x_ix_j} + \sum_{i=1}^n\tilde{b}^{\eps}_i\varphi^{\eps}_{x_i}
+ c^{\eps}\varphi^{\eps} = 0, \label{adj}
\end{equation}
or equivalently,
\begin{equation}
\frac{1}{2}\sum_{i,j=1}^na_{ij}\varphi^{\eps}_{x_ix_j} + \sum_{i=1}^n\frac{\tilde{b}^{\eps}_i}{\eps^2}\varphi^{\eps}_{x_i}
+ \frac{c^{\eps}}{\eps^2}\varphi^{\eps} = 0, \label{adj1}
\end{equation}
where,
\begin{eqnarray*}
\tilde{b}_i^{\eps} &\defn& - b_i + \eps^2\sum_{j=1}^n(a_{ij})_{x_j},\ \mbox{for}\ i=1,2,\cdots,n,
 \\
c^{\eps} &\defn& \frac{\eps^2}{2}\sum_{i,j=1}^n(a_{ij})_{x_ix_j} - \sum_{i=1}^n(b_i)_{x_i}.
\end{eqnarray*}
Define
\begin{displaymath}
A^{\eps} \defn \frac{\eps^2}{2}\mbox{tr}(a\nabla^2) + \langle\tilde{b}^{\eps}, \nabla\rangle.
\end{displaymath}
Letting $W^{\eps}(x) \defn -\eps^2\mbox{ln}(\varphi^{\eps}(x))$, $W^{\eps}$ is  seen
to satisfy
\begin{equation}
A^{\eps}W^{\eps}(x) - \frac{1}{2}DW^{\eps}(x)^Ta(x)DW^{\eps}(x) -
\eps^2c^{\eps}(x) = 0, \label{HJB}
\end{equation}
or, equivalently,
\begin{equation}
\frac{\eps^2}{2}\sum_{i,j=1}^na_{ij}W^{\eps}_{x_ix_j} + \min_{u \in
\Rd}\Big[(\tilde{b}^{\eps} - u)^TDW^{\eps} + \frac{1}{2}u^Ta^{-1}u -
\eps^2c^{\eps}(x)\Big] = 0. \label{HJB1}
\end{equation}
Observe that this is the HJB equation for a certain ergodic
control problem.

\spa

\noi \bo{Lemma 1:}  $\bar{W}^{\eps}(\cdot) \defn W^{\eps}(\cdot) - W^{\eps}(0), \eps \in (0, 1)$, is relatively compact
in $C(\Rd)$.

\spa

\noi \bo{Proof:} In view of (\ref{adj}) and the Harnack estimates of
Theorem $5.2$ of \cite{Met}, we have

\begin{equation}
\|D\bar{W}^{\eps}(x)\| = \|DW^{\eps}(x)\| \leq C \ \forall x\in\Rd \ \label{harn}
\end{equation}
for a constant $C$ depending only on $\lambda,
\|a\|_{C^4_b(\Rd)}, \|b\|_{C^2(\Rd)}$. Since $\bar{W}^{\eps}(0) =
0$, the claim follows by the Arzela-Ascoli theorem. \hfill $\Box$

\spa

Now note that the minimum over $u$ in (\ref{HJB1}) is attained at $u
= a(x)DW^{\eps}(x)$, which is bounded by (\ref{harn}). Thus without
loss of generality, for purposes of analysis  we may a priori
restrict this minimization to a closed bounded set $\Gamma$. Now
letting $\eps\downarrow 0$, standard arguments from the theory of
viscosity solutions (\cite{CrLi0}, Prop. VI.1) tell us that along an appropriate
subsequence, $\bar{W}^{\eps} \rightarrow W$ uniformly on compacts, where $W$
is  a viscosity solution to the p.d.e.
\begin{equation}
\min_{u\in\Gamma}\Big[(- b(x) - u)^TDW +
\frac{1}{2}u^Ta^{-1}u\Big] = 0. \label{HJB2}
\end{equation}
In view of (\ref{harn}), $W$ is in fact Lipschitz. For the
deterministic control system
\begin{equation}
\dot{x} = - b(x) - u(t), \label{detcon}
\end{equation}
define
\begin{displaymath}
\tilde{W}(x,t) =
\inf_{x(0)=x}\Big[\frac{1}{2}\int_0^tu(s)^Ta^{-1}(x(s))u(s)ds +
W(x(t))\Big]
\end{displaymath}
where the infimum is over all measurable and locally
square-integrable $u(\cdot)$. Since $W$ is Lipschitz, $\tilde{W}$
will be a Lipschitz continuous viscosity solution \cite{CrLi0} of
\begin{eqnarray*}
\tilde{W}_t(x,t) &=& \inf_{u \in \Gamma}\Big[(-b(x)
-u)^TD\tilde{W}(x,t) + \frac{1}{2}u^Ta^{-1}(x)u\Big], \\
&& \ \ \ \ \ \ \ \ \ \ \ \ \ \ 0 < t < T, \ \tilde{W}(x, 0) = W(x).
\end{eqnarray*}
From Theorem VI.1 of \cite{Cra}, it follows that this equation has
a unique viscosity solution and thus $W = \tilde{W}$. That is, $W$
is a stationary solution to the above p.d.e, which is also its
unique viscosity solution. In particular, $W$ satisfies:
\begin{equation}
W(x) = \inf_{x(0)=x}\Big[\frac{1}{2}\int_0^tu(s)^Ta^{-1}(x(s))u(s)ds +
W(x(t))\Big], \label{stat-hjb}
\end{equation}
where the infimum is over all $(x(\cdot), u(\cdot))$ satisfying
(\ref{detcon}) with $u(\cdot)$ locally square-integrable. We use
this to establish several additional properties of $W$ in the next
section, leading to our main result.

\spa

\section{Proof of Theorem 2}

We  proceed through a sequence of lemmas.

\spa

\noi \bo{Lemma 2:} $W(x) \geq 0$.

\spa

\noi \bo{Proof:} Suppose $\bar{W}^{\eps_n} \rightarrow W$ uniformly on
compacts and suppose there exists $x_0 \neq 0$ such that $W(x_0) =
-\delta$ for some $\delta
> 0$. Let $B_r(x)$ denote the open ball of radius $r$ centered at $x$.
Then for a sufficiently small $r$, we have
\begin{displaymath}
\bar{W}^{\eps_n}(x) < -\frac{\delta}{2} \ \mbox{for} \ x \in
B_r(x_0); \ \bar{W}^{\eps_n}(x) > -\frac{\delta}{4} \ \mbox{for} \
x \in B_r(0).
\end{displaymath}
So for $n$ sufficiently large,
\begin{eqnarray*}
\frac{\mu^{\eps_n}(B_r(0))}{\mu^{\eps_n}(B_r(x_0))} &=&
\frac{\int_{B_{r}(0)}\varphi^{\eps_n}(x)dx}{\int_{B_r(x_0)}\varphi^{\eps_n}(x)dx}
\\
&=&
\frac{\int_{B_r(0)}e^{-\frac{W^{\eps_n}(x)}{\eps_n^2}}dx}{\int_{B_r(x_0)}e^{-\frac{W^{\eps_n}(x)}{\eps_n^2}}dx}
\\
&=&
\frac{\int_{B_r(0)}e^{-\frac{\bar{W}^{\eps_n}(x)}{\eps_n^2}}dx}{\int_{B_r(x_0)}e^{-\frac{\bar{W}^{\eps_n}(x)}{\eps_n^2}}dx}
\\
&\leq&
\frac{|B_r(0)|e^{\frac{\delta}{4\eps_n^2}}}{|B_r(x_0)|e^{\frac{\delta}{2\eps_n^2}}}
\rightarrow 0 \ \mbox{as} \ n\uparrow\infty.
\end{eqnarray*}
This contradicts Theorem 1, proving the claim. \hfill $\Box$

\spa

\noi \bo{Lemma 3:} For $\{\eps_n\}$ as above,
$\lim_{n\uparrow\infty}W^{\eps_n}(0) = 0$.

\spa

\noi \bo{Proof:} We may write
\begin{displaymath}
\varphi^{\eps_n}(x) =
\varphi^{\eps_n}(0)h_{\eps_n}(x)e^{-\frac{W(x)}{\eps_n^2}}
\end{displaymath}
where $\eps_n^2\mbox{ln}(h_{\eps_n}(x)) \rightarrow 0$ uniformly on compact
sets. Let $A \subset \Rd$ denote a compact set with $0$ in its
interior. Then
\begin{displaymath}
\eps_n^2\mbox{ln}\Big(\int_A\varphi^{\eps_n}(x)dx\Big) =
\eps_n^2\mbox{ln}\Big(\int_A\varphi^{\eps_n}(0)h_{\eps_n}(x)e^{-\frac{W(x)}{\eps_n^2}}dx\Big).
\end{displaymath}
We have
\begin{displaymath}
\left(\int_Ae^{-\frac{W(x)}{\eps_n^2}}dx\right)^{\eps_n^2}
\stackrel{n\uparrow\infty}{\rightarrow} \mbox{ess.sup}\{e^{-W(x)}: x \in
A\}
\end{displaymath}
and
\begin{eqnarray*}
\lefteqn{ \eps_n^2\mbox{ln}\Big(\varphi^{\eps_n}(0)\Big) +
\eps_n^2\mbox{ln}\Big(\int_Ae^{-\frac{W(x)}{\eps_n^2}}dx\Big) +
\eps_n^2\mbox{ln}\Big(\inf_Ah_{\eps_n}(x)\Big) } \\
&\leq&
\eps_n^2\mbox{ln}\Big(\int_A\varphi^{\eps_n}(0)h_{\eps_n}(x)e^{-\frac{W(x)}{\eps_n^2}}dx\Big)
\\
&\leq& \eps_n^2\mbox{ln}\Big(\varphi^{\eps_n}(0)\Big) +
\eps_n^2\mbox{ln}\Big(\int_Ae^{-\frac{W(x)}{\eps_n^2}}dx\Big) +
\eps_n^2\mbox{ln}\Big(\sup_Ah_{\eps_n}(x)\Big).
\end{eqnarray*}
Therefore
\begin{displaymath}
-\eps_n^2\mbox{ln}\Big(\int_A\varphi^{\eps_n}(x)dx\Big) +
\eps_n^2\mbox{ln}\Big(\varphi^{\eps_n}(0)\Big)
\stackrel{n\uparrow\infty}{\rightarrow} \inf\{W(x): x \in A\} = 0.
\end{displaymath}
The first term on the left goes to zero by Theorem $1$. Thus so does
the second, proving the claim. \hfill $\Box$

\spa

\noi \bo{Lemma 4:} $\displaystyle{\lim_{\|x\|\rightarrow\infty}}W(x) = \infty$.

\spa

\noi \bo{Proof:} Let $K > 0$, to be prescribed later and $\delta <
\frac{2\alpha}{\beta}$ for $\alpha, \beta$ as in (\ref{assum}),
$\theta, c > 0, m_1 \geq 2, m_2 \geq 4$.  Define  a $C^2$ function
$U^{\eps_n}$ on $\Rd$ satisfying:
\begin{eqnarray*}
U^{\eps_n}(x) &=& e^{\frac{\delta\|x\|^{\beta}}{\eps_n^2}} \ \ \ \
\ \ \ \ \ \ \ \ \ \ \forall \ \|x\|
> K, \\
&=& 1 \ \ \ \ \ \ \ \ \ \ \ \ \ \ \ \ \ \ \ \ \forall \ \|x\| \leq K
-
\theta\eps_n^2, \\
|\LA^{\eps_n}U^{\eps_n}(x)| &<&
\frac{c}{(\theta\eps_n^2)^{m_1}}e^{\frac{m_2\delta
K^{\beta}}{\eps_n^2}} \ \ \ \forall \ \|x\| \leq K.
\end{eqnarray*}
Then for $\|x\| \geq K$,
\begin{eqnarray*}
\LA^{\eps_n}U^{\eps_n}(x) &=& \delta\beta\|x\|^{\beta -
1}e^{\frac{\delta\|x\|^{\beta}}{\eps_n^2}}\Big[\frac{\sum_{i=1}a_{ii}}{2\|x\|}
+ \frac{\beta - 1}{2\|x\|^3}\sum_{i,j=1}^na_{ij}x_ix_j \\
&& + \ \frac{1}{\eps_n^2}\left(\frac{\delta\beta}{2}\|x\|^{\beta -
3}\sum_{i,j=1}^na_{ij}x_ix_j + b(x)^T\frac{x}{\|x\|}\right)\Big].
\end{eqnarray*}
In view of (\ref{assum}), we may choose $K$ large enough such that
the r.h.s.\ above is $< 0$. Then by Proposition 2.3 of \cite{Met},
$\int |\LA^{\eps_n}U^{\eps_n}|d\mu^{\eps_n} < \infty$ and

\begin{displaymath}
\int_{\Rd} |\LA^{\eps_n}U^{\eps_n}|d\mu^{\eps_n} \leq 2\int_{B_K(0)}
|\LA^{\eps_n}U^{\eps_n}|d\mu^{\eps_n} \leq
\frac{2c}{(\theta\eps_n^2)^{m_1}}e^{\frac{m_2\delta
K^{\beta}}{\eps_n^2}}.
\end{displaymath}
Choose $K_1 > K$ such that
\begin{displaymath}
|\LA^{\eps_n}U^{\eps_n}| \geq e^{\frac{\delta
\|x\|^{\beta}}{2\eps_n^2}} \ \ \ \ \ \ \forall \ \|x\| \geq K_1.
\end{displaymath}
Then \begin{displaymath} \int_{B_{R}^c(0)}e^{\frac{\delta
\|x\|^{\beta}}{2\eps_n^2}}d\mu^{\eps_n} \leq
\frac{2c}{(\theta\eps^2)^{m_1}}e^{\frac{m_2\delta
K^{\beta}}{\eps_n^2}}, \ \ R > K_1..
\end{displaymath}
Hence for $R > K_1$,
\begin{displaymath}
\mu^{\eps_n}(\bar{B}^c_R(0)) \leq
\frac{2c}{(\theta\eps_n^2)^{m_1}}e^{\frac{2m_2\delta K^{\beta} -
\delta R^{\beta}}{2\eps_n^2}}
\end{displaymath}
Let $\bar{x} \in \bar{B}^c_R(0)$ and $\Theta$ an open and bounded neighborhood
thereof. Then
\begin{eqnarray*}
-\eps_n^2\mbox{ln}([\inf_{x\in\Theta}\varphi^{\eps_n}(x)]|\Theta|)
&\geq& -\eps_n^2\mbox{ln}(\int_{\bar{B}_R(0)}\varphi^{\eps_n}(x)dx)
\\
&\geq& -\eps^2_n\mbox{ln}(2c) + m_1\eps_n^2\mbox{ln}(\theta\eps_n^2) +
\frac{\delta R^{\beta}}{2} - m_2\delta K^{\beta},
\end{eqnarray*}
leading to, in view of Lemma $3$,
\begin{displaymath}
\sup_{\Theta}W(x) \geq \frac{\delta R^{\beta}}{2} - m_2\delta
K^{\beta}.
\end{displaymath}
For $R$ large enough and $\Theta$ small enough, we have by the
Lipschitz continuity of $W$ that $W(\bar{x}) \geq \frac{\delta
R^{\beta}}{4}$, which implies the result in view of our arbitrary
choice of $\bar{x}$. \hfill $\Box$

\spa

\noi \bo{Lemma 5:}\ $W(x) > 0$ for $x \neq 0$.

\spa

\noi \bo{Proof:} Suppose $W(x) = 0$ for some $x \neq 0$.
Considering $x(0) = x$ in (\ref{detcon}), we have
\begin{displaymath}
0 = W(x) = \inf\Big[\frac{1}{2}\int_0^Tu(t)^Ta(x(t))^{-1}u(t)dt +
W(x(T))\Big] \ \forall \ T > 0,
\end{displaymath}
implying that the infimum is in fact attained at $u(\cdot) \equiv
0$ and  the corresponding $x(\cdot)$ satisfies $W(x(T)) = 0 \
\forall \ T > 0$. But in view of (\ref{detcon}) with $u(\cdot)
\equiv 0$ and our hypotheses on $b(\cdot)$, $x(T)\uparrow\infty$, leading
to $W(x(T))\uparrow\infty$, a contradiction. The result follows.
\hfill $\Box$

\spa

\noi Thus we have:

\spa

\noi \bo{Lemma 6:}\ $W(x) =
\inf\Big[\frac{1}{2}\int_0^{\infty}u(t)^Ta^{-1}(x(t))u(t)dt\Big]$
where the infimum is over all $(x(\cdot), u(\cdot))$ satisfying
(\ref{detcon}) with the additional restrictions: $x(0) = x$ and $x(t)
\stackrel{t\uparrow\infty}{\rightarrow} 0$.

\spa

\noi \bo{Proof:} We topologize the space of locally square-integrable
$u(\cdot)$ as follows: For $T > 0$, let $L^2_w[0,T]$ denote the
space $L^2([0,T]; \Rd)$ with weak$^*$ topology. Equip the set of
admissible $u(\cdot)$ with the coarsest topology that renders
continuous the map $u(\cdot) \rightarrow u(\cdot)|_{[0,T]} \in L^2_w[0,T]$
for all $T > 0$. Now consider the minimization problem of minimizing
over $u(\cdot) \in L^2_w[0,T]$ the functional
\begin{equation}
\frac{1}{2}\int_0^Tu(t)^Ta(x(t))^{-1}u(t)dt + W(x(T)), \label{lemsix}
\end{equation}
where $x(\cdot), u(\cdot)$ are related through (\ref{detcon}),
with $x(0) = x$ (say).
Let
\begin{displaymath}
 \U_T=\lbrace u\in L^2_w[0,T]\ |\ u(t)\in\Gamma\ \mbox{a.e.}\rbrace
\end{displaymath}
It is easy to verify that (\ref{lemsix}) is a lower
semicontinuous functional which will attain its minimum over a
compact set $\U_T$ of $L^2_w[0,T]$. Let $T'
> T$. Then a standard dynamic programming argument shows that
\begin{displaymath}
u(\cdot) \in \U_{T'} \Longrightarrow u(\cdot)|_{[0,T]} \in \U_T.
\end{displaymath}
Define $\U^*_T \defn$ the set of $u(\cdot) \in \U(\equiv\U_{\infty})$ such that
$u(\cdot)|_{[0,T]} \in \U_T$. Then as $T\uparrow\infty$, it is a
family of nested decreasing compact subsets of $\U$ and therefore
has a non-empty intersection. Take $u^*(\cdot)$ in this intersection
and let $x^*(\cdot)$ denote the corresponding trajectory of
(\ref{detcon}). Then
\begin{displaymath}
W(x) = \frac{1}{2}\int_0^Tu^*(t)a(x^*(t))^{-1}u^*(t)dt + W(x^*(T)) \
\forall \ T,
\end{displaymath}
implying that $W(x^*(T)) \downarrow W^*$ for some $W^* \geq 0$.
Suppose $W^* > 0$. Then there exist $0 < \delta < R$ such that
$\delta < x^*(s) \leq R \ \forall \ s$. By Lemma 3.1 of
\cite{Fle}, we then have
\begin{displaymath}
W(x) \geq \frac{1}{2}\int_0^Tu^*(t)a^{-1}(x^*(t))u^*(t)dt \rightarrow
\infty,
\end{displaymath}
a contradiction. Hence $W^* = 0$. By Lemma $5$, we have $x^*(T)
\stackrel{T\uparrow\infty}{\rightarrow} 0$. The rest is easy. \hfill
$\Box$

\spa

\noi \bo{Proof of Theorem 2:} Immediate from the above lemmas.
\hfill $\Box$

\spa

Let $R_\eps=\eps^de^{W/\eps^2}\varphi^\eps$. In \cite{Day1} it was shown that for a certain class of drifts, $W$ is $C^1$ in an open, dense, connnected set $\bo{G}\subset\Rd$ and
\begin{center}
 $R_\eps\rightarrow R_0$ \ \ \ \ \ \ uniformly on compact subset in $\bo{G}$,
\end{center}
where $R_0$ satisfies the equation
\begin{equation}
 <b+a\nabla W,\nabla R_0> + (div(W)+\frac{1}{2}\bigtriangleup W)R_0=0\label{dayreg}
\end{equation}
in $\bo{G}$.
The same result can be generalized for the class of drifts considered here, which
subsume the drifts considered in \cite{Day1} for $\beta<1$. This is because the asymptotic stability of $0$ for (\ref{ode}) implies that the Jacobian matrix $Db(0)$ will have eigenvalues with strictly negative real parts (see, e.g., section 1.3 of \cite{Guk}). We state it as theorem
below and omit the proof which is identical to that of Theorem $3$ \cite{Day1}. (Of course, the equation (\ref{dayreg}) gets suitably modified.)

\spa

\noi \bo{Theorem 3:} There exists a positive function $R_0\in C^1(\bo{G})$
such that $R_\eps\stackrel{\eps\downarrow 0}{\to}R_0$.

\section{Extension to multiple equilibria}

Now consider the case when (\ref{ode}) has finitely many equilibria
$x_1, \cdots, x_J$ (say), and no other $\omega$-limit sets. In this
case, one can mimic the arguments of Theorem 1 to claim that all
limit points of $\mu^{\eps}$ as $\eps\downarrow 0$
concentrate on the set $\{x_1, \cdots, x_J\}$. Furthermore, we again
have $W$ as the Lipschitz continuous subsequential limit of
$-\eps^2\mbox{ln}\Big(\frac{\varphi^{\eps}(x)}{\varphi^{\eps}(0)}\Big)$
as $\eps\downarrow 0$, which satisfies (\ref{stat-hjb}). Also, by
arguments of Lemma 4, we have $\lim_{\|x\|\uparrow\infty}W(x) =
\infty$. The proof of Lemma 2  can be adapted to show that the limit
points of $\mu^{\eps}$ as $\eps\downarrow 0$ will in fact concentrate
on the set Argmin$(W(x): x \in \{x_1, \cdots, x_J\})$.

Let $(\hat{x}(\cdot), \hat{u}(\cdot))$ denote an optimal pair with
$\hat{x}(0) = x$, as in Lemma 6. Then for $t, T > 0$,
\begin{equation}
W(\hat{x}(t)) =
\frac{1}{2}\int_t^{t+T}\hat{u}(s)^Ta^{-1}(\hat{x}(s))\hat{u}(s)ds +
W(\hat{x}(t + T)). \label{DPcont}
\end{equation}
Thus $W(\hat{x}(t))$ is non-increasing with $t$, implying in
particular that $\hat{x}(\cdot)$ remains bounded. Thus
$\int_t^{t+T}\hat{u}(s)^Ta^{-1}(\hat{x}(s))\hat{u}(s)ds$ is bounded.
By Banach-Alaoglu theorem,  $\hat{u}(t + \cdot)|_{[0,T]}, t \geq 0,$
is relatively compact in $L^2_w[0, T]$. By a standard argument using
Arzela-Ascoli theorem, one verifies that $\hat{x}(t +
\cdot)|_{[0,T]}$ is also relatively compact in $C([0,T]; \Rd)$. Let
$(x^*(\cdot), u^*(\cdot))$ denote a subsequential limit of
$(\hat{x}(t + \cdot), \hat{u}(t + \cdot))$ in $C([0,T]; \Rd)\times
L^2_w[0,T]$ as $t\uparrow\infty$. Then by the lower semicontinuity
of the map
\begin{displaymath}
\hat{u}(\cdot) \in L^w_2[0,T] \rightarrow
\int_0^T\hat{u}(s)^Ta^{-1}(\hat{x}(s))\hat{u}(s)ds,
\end{displaymath}
we have
\begin{equation}
W(x^*(0)) \geq \frac{1}{2}\int_0^{T}u^*(s)^Ta^{-1}(x^*(s))u^*(s)ds +
W(x^*(T)). \label{DPcont2}
\end{equation}
But $W(\hat{x}(t))$ is monotonically decreasing by virtue of
(\ref{DPcont}), whence $W(x^*(0)) = W(x^*(T)) = c \defn
\lim_{t\uparrow\infty}W(\hat{x}(t))$. Thus $u^*(\cdot) \equiv 0$
a.e. We drop the `a.e.' without loss of generality. It follows that
$\hat{x}(t + \cdot)$ converges as $t\uparrow\infty$ to
$H=\cap_{t>0}\overline{\{\hat{x}(t+s),s\geq 0\}}$, which is a positively invariant set of the o.d.e.
\begin{displaymath}
\dot{x}(t) = -b(x(t)).
\end{displaymath}
 A similar argument can be given to
establish negative invariance of this set as well, implying that it
converges to an invariant set of the above o.d.e. But only such sets
are the equilibria of (\ref{ode}). Since $H$ is compact connected
(being intersection of such), it must converge to a single
equilibrium. Letting $t = 0$ and $T\uparrow\infty$ in
(\ref{DPcont}), we then have
\begin{displaymath}
W(x) =
\frac{1}{2}\int_0^{\infty}\hat{u}(t)^Ta^{-1}(\hat{x}(t))\hat{u}(t)dt
+ W(x_i)
\end{displaymath}
for some $1 \leq i \leq J$. Since we also have
\begin{displaymath}
W(x) = \inf\Big[\frac{1}{2}\int_0^Tu(t)^Ta^{-1}(x(t))u(t)dt +
W(x(T))\Big],
\end{displaymath}
it follows by a straightforward argument that:

\spa

\noi \bo{Theorem 4:}
\begin{equation}
W(x) = \min_{1\leq i\leq J}\inf_{(x(\cdot),u(\cdot)):x(t)\to x_i}
\Big[\frac{1}{2}\int_0^{\infty}u(t)^Ta^{-1}(x(t))u(t)dt +
W(x_i)\Big]. \label{star}
\end{equation}

\spa

In particular, it follows that $W$ attains its minimum at one or
more of the $x_i$'s. As $W$ needs be specified only up to an additive
factor, we may assume without loss of generality that its minimum
value is zero.

Note, however, that unlike in the single equilibrium case, the
uniqueness of $W$, obtained as a \textit{subsequential} limit, is
not immediate.

\spa

\noi \textbf{Theorem 5:} $W$ is uniquely specified as the Lipschitz
function satisfying (\ref{star}) and the condition $\min W = 0$.

\spa

\noi \textbf{Proof:} Note that by (\ref{star}), we have by arguments
similar to those used in Lemma 6 that for each $x$, there is an
optimal pair $(x^*(\cdot), u^*(\cdot))$ such that $x^*(t) \rightarrow x_{i^*}$
for some $i^*$ that will depend on $x$, and
\begin{displaymath}
W(x) = \frac{1}{2}\int_0^{\infty}u^*(t)^Ta^{-1}(x^*(t))u^*(t)dt +
W(x_{i^*}).
\end{displaymath}
If there is more than one $i^*$ for which this holds, we choose one according to some pre-specified rule. In this case, write $x \longrightarrow x_{i^*}$. We may also have
$x_i  \longrightarrow x_j$ for some $i \neq j$, with $W(x_j) <
W(x_i)$. In this case, write $x_i \Longrightarrow x_j$. If we draw a
directed graph with nodes $\{x_1, \cdots, x_J\}$ with a directed
edge from $x_i$ to $x_j$ whenever $x_i \Longrightarrow x_j$, we get
a forest of rooted trees, say, $T_1, \cdots, T_K$. For each $T_i$,
let $\hat{x}_i$ denote its root. Let $O_i
\defn \{x: x \longrightarrow x_j \ \mbox{for some} \ x_j \in T_i\}$. Then
$\Rd$ is the disjoint union of the $O_i$'s. On each $O_i$, $W$ is
completely specified in terms of $W(\hat{x}_i)$ by successive
application of (\ref{star}) as follows: First use (\ref{star}) to
obtain values of $W$ at nodes $\{x_j\}$ in $T_i$ one removed from
$\hat{x}_i$, then repeat the same for nodes two removed, and so on
till $W$ is defined for all nodes in $T_i$. Then define it for $x
\in O_i$ by using (\ref{star}) again. In particular, $W$ is
completely specified for those $T_i$ for which $W(\hat{x}_i) = 0$.
For others, it is in principle specified up to an additive scalar,
because the value of $W(\hat{x}_i)$ is not specified. But we have
the additional restriction that $W$ be continuous (in fact,
Lipschitz) over the whole of $\Rd$, whence this choice is also
unique. Uniqueness of $W$ as $\lim_{\eps\downarrow
0}\Big(-\eps^2\mbox{ln}(\varphi^{\eps}(\cdot))\Big)$ follows. \hfill
$\Box$
\spa

\noindent \textbf{Remark:} The problem of finding stationary density is one of finding the principal eigenfunction of an operator. A counterpart of the above in this general framework, albeit in a discrete set-up, appears in \cite{Sheu2}.\\


\begin{thebibliography}{99}

\bibitem{Ari1} ARISAWA, M.\ (1997) ``Ergodic problem for the
Hamilton-Jacobi-Bellman equation. I. Existence of the ergodic attractor",
\textit{Ann.\ Inst.\ Henri Poincar\'{e}} 14, 415-438.

\bibitem{Ari2} ARISAWA, M.\ (1998) ``Ergodic problem for the Hamilton-Jacobi-Bellman
equation. II", \textit{Ann.\ Inst.\ Henri Poincar\'{e}} 15, 1-24.


\bibitem{CrLi0} CRANDALL, M.\ G.; LIONS, P.\ L.\ (1983)
``Viscosity solutions of Hamilton-Jacobi equations", \textit{Trans.\
American Math.\ Soc.\ } 277, 1-38.

\bibitem{Cra} CRANDALL, M.\ G.; LIONS, P.\ L.\ (1987) ``Remarks on the existence
and uniqueness of unbounded viscosity solutions of Hamilton-Jacobi equations",
\textit{Illinois J.\ Mathematics} 31, 665-688.

\bibitem{DaDa} DAY, M.\ V.; DARDEN, A.\ T.\ (1985) ``Some regularity results on the Ventcell-Freidlin
quasi-potential function``, \textit{Appl.\ Math.\ Optim.} 13, 259-282.

\bibitem{Day1} DAY, M.\ V.\ (1987), ``Recent progress on the small parameter exit promlem",
\textit{Stochastics} 20, 121-150.


\bibitem{EcRu} ECKMANN, J.\ P.; RUELLE, D.\ (1985) ``Ergodic
theory of chaos and strange attractors", \textit{Reviews of Modern
Physics} 57, 617-656.

\bibitem{Eth} ETHIER, S.\ E.; KURTZ, T.\ G.\ (1986)
\textit{Markov Processes: Characterization and Convergence}, John
Wiley, New York.

\bibitem{Fle} FLEMING, W.\ H.\ (1978) ``Exit probabilities and optimal stochastic control",
\textit{Appl.\ Math.\ Optim.\ } 4, 329-346.

\bibitem{FoYo} FOSTER, D.; YOUNG, H.\ P.\ (1990), ``Stochastic
evolutionary game dynamics", \textit{Theoretical Population Biology}
38, 219-232.


\bibitem{Fre} FREIDLIN,  M.\ I.; WENTZELL,  A.\ D.\ (1998) ``Random Perturbations
of Dynamical Systems",  Springer Verlag,  New York.

\bibitem{Frie} FRIEDMAN,  A.\  (1987) ``Stochastic Differential Equations``, Vol 2,
\textit{Academic Press}, New York.


\bibitem{GihSko} GIHMAN, I.\ I.; SKOROKHOD, A.\ V.\ (1972) \textit{Stochastic Differential Equations}, Spriinger Verlag, New York.

\bibitem{GiTu} GILBARG, D.; TRUDINGER, N.\ S.\ (1977) ``Elliptic Partial differentail equation of second order",\textit{Springer-Verlag},Berlin.

\bibitem{Guk} GUCKENHEIMER, J.; HOLMES, P.\ (1983) \textit{Nonlinear Oscillations, DDynamical Systems, and Bifurcations of Vector Fields}, Springer Verlag, New York.

\bibitem{Khas} KHASMINSKII, R.\ Z.\ (1960) ``Ergodic properties of recurrent diffusion processes and stabilization of the soluution to the Cauchy problrm of parabolic equations", \textit{Theory of Probability and Appl.} 2, 179-196.

\bibitem{Khas2} KHASMINSKII, R.\ Z.\ (1980) \textit{Stochastic Stability of Differential Equations}, Sijthoff and Noordhoff, Leyden, The Netherlands.


\bibitem{Met} METAFUNE, G.; PALLARA, D.; RHANDI, A.\ (2005) ``Global properties of invariant measures",
\textit{J.\ Functional Analysis} 223, 396-424.

\bibitem{Mik} MIKAMI, T.\ ``Asymptotic expansions for the invariant density of a Markov process with a small parameter",
\textit{Ann.\ Inst.\ Henri Poincar\'{e}}  24, 403-424.

\bibitem{Sas} SASTRY, S.\ S.\ (1983) ``The effects of small noise on
implicitly defined nonlinear dynamical systems", \textit{IEEE
Trans.\ on Circuits and Systems} CAS-30, 651-663.

\bibitem{She} SHEU,  S.--J.\  (1986) ``Asymptotic behavior of the invariant density of a diffusion
Markov process with small diffusion",  \textit{SIAM J. \ Math. \ Analysis} 17,  451-460.

\bibitem{Sheu2} SHEU, S.--J.; AND WENTZELL, A.\ D.\ (1999) ``On the solutions of the equation arising from the singular limit of some eigen problems", in `\textit{Stochastic Analysis, Control, Optimization and Applications}' (W.\ M.\ McEneaney, G.\ George Yin, Q.\ Zhang, eds.) Birkh\"{a}user, Boston, 135-150.

\bibitem{stet} STETTNER, L.\ (1989) ``Large deviations of invariant measures for de-
generate diffusions",\textit{Probability and Mathematical Statistics} 10,
 93-105.

\bibitem{Yosi} YOSIDA, K.\ (1980) \textit{Functional Analysis} (6th ed.), Springer Verlag, Berlin-Heidelberg.
\end{thebibliography}
\end{document}


\maketitle

The uniqueness argument in the proof of Theorem 5, p.\ 483, of \cite{bisbor} is flawed. We give here a corrected proof.\\

\noi Define
$$ S_{0T}(\phi)=\frac{1}{2}\int_0^T(\dot\phi(s)-b(\phi(s)))a^{-1}(\phi(s))(\dot\phi-b(\phi(s)))ds.$$
and
$$V_i(x)=\inf\{S_{0T}(\phi): \ \phi(0)=x_i, \ \phi(T)=x\},$$
where the quantities are defined for absolutely continuous paths.
\begin{lemma}
 Under the stated assumption in \cite{bisbor}, $V_i(x)\to\infty$ as $\|x\|\to\infty$ for $i\leq i\leq J$.
\end{lemma}

\noi{\bf Proof.}\ Consider the stochastic differential equation
$$dY^\eps(t)= b(Y^\eps(t))dt+\eps dB(t),$$
where $B(\cdot)$ is a $d$-dimensional Brownian motion. Then following the arguments in \cite{bisbor}, page 477--479, we have a Lipschitz continuous viscosity solution $U(\cdot)$
to the p.d.e.
$$ \min_{u\in\R^d}[(-b(x)-u).\nabla U +\frac{1}{2}\|u\|^2]=0.$$
Also arguing as in Lemma 2, \cite{bisbor} we have $U(\cdot)\geq 0$. As in \cite{bisbor}, we define
$$\tilde U(x, t)=\inf_{x(0)=x}[\frac{1}{2}\int_0^t\|u(s)\|^2ds + U(x(t))],$$
where the infimum is taken over all locally square integrable $u(\cdot)$ satisfying $\dot{x}(s)=-b(x(s))-u(s),\ x(0)=x$. One can easily check that $\tilde U(x,t)$
is continuous on $\R^d\times[0 ,T]$ and Lipschitz in space variable uniformly w.r.t. $t\in [0, T]$. Also $\tilde U(x,t)$ is a viscosity solution
to the equation
$$\tilde{U}_t(x, t)=\min_{u\in\R^d}[(-b(x)-u).\nabla \tilde U +\frac{1}{2}\|u\|^2].$$
Therefore applying Theorem VII.2. in \cite{crandall} (see  also the Appendix) we have
\begin{eqnarray}
 U(x)=\inf_{x(0)=x}[\frac{1}{2}\int_0^T\|u(s)\|^2ds + U(x(T))],\label{hjb}
\end{eqnarray}
where the infimum is over all $(x(\cdot),u(\cdot))$ satisfying
$$ \dot{x}(t)=-b(x(t))-u(t), \ x(0)=x$$
with $u(\cdot)$ locally square integrable. Furthermore, applying Lemma 4 in \cite{bisbor} we have $U(x)\to\infty$ as $\|x\|\to\infty$. Suppose there exists
a sequence $\{y_n\}$ such that $\|y_n\| \uparrow \infty$ and $V_i(y_n)\leq M$ for all $n\geq 1$ and some $0 < M < \infty$. Therefore there exists $\{\phi_n(\cdot)\}$ such that
$$S_{0T_n}(\phi_n)\leq M+1,\ \phi_n(0)=x_i, \ \phi_n(T_n)=y_n.$$
Define $\psi_n(s)=\phi_n(T_n-s)$ for $s\leq T_n$ and
\[u_n(s)=\left\{\begin{array}{ll}
                 -\dot\psi_n(s)-b(\psi_n(s))\ & \mbox{for}\ s\leq T_n
\\
                        0 & \mbox{otherwise}.
                \end{array}
\right. \]
Therefore $\dot\psi_n(s)=-b(\psi_n(s))-u(s)$. Hence from (\ref{hjb}) we have
\begin{eqnarray*}
 U(y_n) &\leq & \frac{1}{2}\int_0^{T_n}\|u_n(s)\|^2 ds + U(x_i)
\\
&=& \Lambda S_{0T_n}(\phi_n) + U(x_i)\leq \Lambda(M+1) + U(x_i),
\end{eqnarray*}
where $\Lambda$ is the maximal eigenvalue of $a(\cdot)$.
This contradicts the fact that $U(x)\to\infty$ as $\|x\|\to\infty$. This proves the claim.\hfill $\Box$

\spa

Since there are only finitely many equilibria, it is easy to see that $V(x):=\min_{1\leq i\leq J}V_i(x)\to\infty$ as $\|x\|\to\infty$. We introduce
the following notations:
$$ \tilde V(x_i,x_j)=\inf\{S_{0T}(\phi):\ \phi(0)=x_i, \phi(T)=x_j, \phi(s)\in\R^d\setminus\cup_{k\neq i.j}\{x_k\}\}.$$
Let $B_R$ denotes the ball of radius $R$ around $0$. For $R$ large, define
\begin{eqnarray*}
\tilde V_R(x_i,x_j) &=& \inf\{S_{0T}(\phi):\ \phi(0)=x_i, \phi(T)=x_j, \\
&& \ \ \ \ \ \ \ \ \ \ \phi(s) \in \bar{B}_R\setminus\cup_{k\neq i.j}\{x_k\} \ \forall \ s \in [0,T]\}.
\end{eqnarray*}

\spa

\noi Let $\rho_1, \rho_0$ be such that $0<\rho_1<\rho_0<\min_{1\leq i, j\leq J}\|x_i-x_j\|$. Let $\Gamma_i, g_i$ denote the boundaries of
the balls of radius $\rho_0, \rho_1$ respectively around $x_i$. Denote $\Gamma=\cup_{1\leq i\leq J}\Gamma_i$ and $g=\cup_{1\leq i\leq J}g_i$.
Define $\tau_0=0, \sigma_n = \inf\{t\geq\tau_n : X^\eps(t)\in\Gamma\}, \tau_n = \min\{t\geq\sigma_{n-1}: X^\eps(t)\in g\}$, where $X^\eps(\cdot)$ satisfies the
stochastic differential equation  (1) in \cite{bisbor}, i.e.,
$$dX^\eps(t)=b(X^\eps(t))dt+\eps\sigma(X^\eps(t))dB(t). $$
 Define the Markov chain $Z_n := X^\eps(\tau_n)$.

\begin{lemma}
 For any $\gamma>0$ there exists $\rho_0$ small such that for any $\rho_2, 0<\rho_2<\rho_0$, there exists $\rho_1, 0<\rho_1<\rho_2,$ such that
for $\eps$ sufficiently small, for all $x\in g_i$ the one-step transition probabilities of $Z_n$ satisfy the inequalities
$$ \exp\left(-\frac{\tilde V(x_i,x_j)+\gamma}{\eps^2}\right)\leq P\left(x,g_j\right)\leq \exp\left(-\frac{\tilde V(x_i,x_j)-\gamma}{\eps^2}\right).$$
\end{lemma}

\noi{\bf Proof.} We choose $R>0$ large enough such that following holds:
\begin{itemize}
 \item Any solution of $\dot{x}(t)=b(x(t))$ that starts in $\bar{B}_R$, stays in $\bar{B}_R$,
\item $\tilde V_R(x_i,x_j) -\frac{\gamma}{2}\leq \tilde V(x_i, x_j)\leq \tilde V_R(x_i,x_j) $ for all $1\leq i, j\leq J$,
\item $\min_{x\in\partial B_R}V(x)> \tilde V(x_i, x_j)+ 4\gamma$ for all $1\leq i, j\leq J$.
\end{itemize}
Therefore applying Lemma 6.2.1 in \cite{freid} we have for $x\in g_i$,
$$ P(x, g_j)\geq \exp\left(-\frac{\tilde V_R(x_i, x_j)+\gamma/2}{\eps^2}\right)\geq \exp\left(-\frac{\tilde V(x_i, x_j)+\gamma}{\eps^2}\right)$$
for $\eps$ small. To get the upper bound we follow the calculations in (\cite{freid}, page 174). Following the calculations on page 174 (estimate (2.6))
of \cite{freid} it is enough to bound $P_x(\tau_1>T)$ by $\exp\left(-\frac{\tilde V(x_i,x_j)-\gamma/2}{\eps^2}\right)$ for a chosen $T>0$ and $x\in\Gamma_i$. At this point
we define $\xi=\min\{t>0 : X^\eps(t)\in \partial B_R\}$. Let $V_0=\min_{x\in\partial B_R} V(x)$. Now for $x\in\Gamma_i$
$$ P_x(\tau_1>T)\leq P_x(\tau_1>T, \xi>T) + P_x(\xi\leq T).$$
Using Lemma 6.1.9 in \cite{freid}, we can choose $T$ large such that $\sup_{x\in\Gamma_i}P_x(\tau_1>T, \xi>T)\leq \exp\left(-\frac{V_0}{\eps^2}\right)$ for $\eps$ small enough.
We fix this choice of $T$.
To estimate the other term we define $\Phi(x)=\{\phi : S_{0T}(\phi)\leq V_0-\gamma/4\}$ for $x\in\Gamma_i$. Now any path starting at $x$ that leaves
$B_R$ by time $T$ must be at least $\eta$ distance away in $C([0,T],\R^d)$ from $\Phi(x)$ for $\rho_0, \eta$ small and $x\in\Gamma_i$. Therefore using the Freidlin-Wentzell large deviations principle
we have (see also \cite{freid}, page 174)
$$\sup_{x\in \Gamma_i}P_x(\xi\leq T)\leq \sup_{x\in\Gamma_i}P_x(\mbox{dist}_{0T}(X^\eps,\Phi(x))\geq \eta)\leq \exp\left(-\frac{V_0-\gamma/2}{\eps^2}\right)$$
for $\eps$ small where $\mbox{dist}_{0T}(\phi_1,\phi_2)=\sup_{[0, T]}\|\phi_1(s)-\phi_2(s)\|$ is the uniform metric on $C([0,T],\R^d)$. This completes the proof.\hfill $\Box$

\spa

\noi Taking $\V=\|x\|^2$, for $\eps$ small we have from \cite{bisbor} (estimate (4) of \cite{bisbor})
\begin{eqnarray}
\mathcal{L}^\eps\V\leq -\delta_2, \ \|x\|\geq K.\label{lyapu}
\end{eqnarray}
for some positive $\delta_2$. Let $D_1=B_K$ and $D=B_{K+1}$.
We need to estimate the hitting time $\tau=\min\{t : X^\eps(t)\in \cup_i\partial g_i\}$. Let $O_i$ denote the closed ball of radius $\rho_1$ around $x_i$.

\begin{lemma}
Let $\KK=\overline{D_1\setminus\cup O_i}$.
Then there exists $M_1>0$ such that $$\sup_{x\in\KK}E_x[\tau]\leq M_1$$ for all $\eps\in(0, \eps_0)$.
\end{lemma}

\noi{\bf Proof.}\ Define $\tau_0=\min\{t : X^\eps(t)\notin D\setminus\cup_i O_i\}$, $\sigma_i=\min\{t\geq\tau_{i-1} : X^\eps(t)\in\partial D_1\}$,
$\tau_i=\min\{t\geq\sigma_i : X^\eps(t)\notin D\setminus\cup_i O_i\}.$ Using (\ref{lyapu}) we can show that
$$ \sup_{x\in\partial D} E_x[\sigma_1]\leq \sup_{x\in\partial D}\frac{\V(x)}{\delta_2}=M_2.$$
Also, Lemma 1.9, pp. 168, of \cite{freid} shows that
$$ \sup_{\KK}E_x[\tau_0]\leq M_3,$$
for $\eps>0$ small. If for some $i$, $X^\eps(\tau_i)\in\cup_i\partial g_i$ then we define $\sigma_{i+1}=\tau_{i+1}=\cdots=\tau_i$.
 By positive recurrence we know that $P(\tau<\infty)=1$. Hence
for $x\in\KK$
\begin{eqnarray*}
 E_x[\tau] &=& E_x[\sum_{i=0}^\infty \tau_i\I_{\{\tau=\tau_i, \tau>\tau_{i-1}\}}]
\\
&=& E_x[\sum_{i=0}^\infty\sum_{k=0}^i(\tau_k-\tau_{k-1})\I_{\{\tau=\tau_i, \tau>\tau_{i-1}\}}] \ [\mbox{with}\ \tau_{-1}=0]
\\
&=& E_x[\sum_{k=0}^\infty\sum_{i=k}^\infty(\tau_k-\tau_{k-1})\I_{\{\tau=\tau_i, \tau>\tau_{i-1}\}}]
\\
&=& E_x[\sum_{k=0}^\infty(\tau_k-\tau_{k-1})\I_{\{ \tau>\tau_{k-1}\}}]
\\
&\leq & \sup_{x\in\KK}E_x[\tau_0] + E_x[\sum_{k=1}^\infty(\tau_k-\tau_{k-1})\I_{\{\tau>\tau_{k-1}\}}]
\\
&=& M_3 + \sum_{k=1}^\infty E_x[(\tau_k-\tau_{k-1})\I_{\{\tau>\tau_{k-1}\}}]
\\
&\leq & M_3 + \sup_{x\in \partial D} E_x[\tau_1] \sum_{k=1}^\infty P(\tau>\tau_{k-1})
\end{eqnarray*}
For $x\in \partial D$ we observe that $E_x[\tau_1]=E_x[\tau_1-\sigma_1+\sigma_1]\leq M_3+M_2$.
Also for $k=1$, $sup_{x\in \KK} P(\tau>\tau_0)<\delta$ for some $\delta<1/2$ and $\eps>0$ small. (This is because then the probability concentrates around the   deterministic trajectories of the o.d.e.\ for $\eps = 0$, which always hit
 $\cup_i\partial g_i$ before hitting $\partial D$.) So for $k\geq 2$,
\begin{eqnarray*}
P(\tau>\tau_{k-1}) &\leq & P(\tau>\sigma_{k-1})\sup_{x\in\partial D_1}P_x(\tau>\tau_0)
\\
&\leq & \delta P(\tau> \tau_{k-2})
\\
&\leq& \delta^{k-1}.
\end{eqnarray*}
Hence combining with the above calculations we have $\sup_{x\in\KK}E_x[\tau]\leq M_1 < \infty$ for all $\eps\in(0, \eps_0)$.\hfill $\Box$

\spa

\noi Let $G(i)$ denote the set of all $i$-graphs  defined on page 177 in \cite{freid}. Define
$$ Z(x_i)=\min_{\chi\in G(i)}\sum_{(m\to n)\in \chi}\tilde V(x_m, x_n).$$

\noi Recall that $\mu^\eps$ is the inavariant measure of $X^\eps(\cdot)$. Proof of the following theorem is same as the proof of Theorem 6.4.1 in \cite{freid}
\begin{theorem}\label{invmeas}
 For any $\gamma>0$ there exists $\rho_1>0$ (can be chosen arbitrary small) such that the $\mu^\eps$-measure of $x_i+B_{\rho_1}$ is
between
$$ \exp\left(-\frac{Z(x_i)-\min_j Z(x_j)\pm \gamma}{\eps^2}\right)$$
for $\eps$ sufficiently small.
\end{theorem}
Now we are ready to prove Theorem 5 in \cite{bisbor}.

\begin{theorem}
The function $W(\cdot)$ obtained in Theorem 4 of \cite{bisbor} is unique and given by
$$ W(x_i)=\min_{1\leq i\leq J}\inf_{(x(\cdot), u(\cdot)): x(t)\to x_i}[\frac{1}{2}\int_0^\infty u(t)^T a^{-1}(x(t))u(t)dt + W(x_i)],$$
where $(x(\cdot), u(\cdot))$ satisfies (16) in \cite{bisbor}.
\end{theorem}

\noi{\bf Proof.}\ The above expression implies that $W(\cdot)$ is uniquely determined by the values $W(x_i), 1\leq i\leq J$. Recall that there
exists a subsequence $\{\eps_n\}$ such that $-\eps_n^2\ln (\varphi^{\eps_n}(x))\to W(x)$ uniformly on compact subsets of $\R^d$ as $\eps_n\to 0$ where $\varphi^\eps(\cdot)$
denotes the invariant density of $\mu^\eps$. To prove the uniqueness it is enough to show that the values of $W(x_i), 1\leq i\leq J,$ are uniquely
determined. Let $\gamma>0$ be small. Choose $\rho_1$ small enough so that $|W(x)-W(x_i)|<\gamma$ for all $\|x-x_i\|\leq\rho_1$.
Therefore applying Theorem \ref{invmeas} we have
$$-\eps_n^2\ln[\sup_{B_{\rho_1(x_i)}}\varphi^{\eps_n}|B_{\rho_1}(x_i)|]\leq Z(x_i)-\min_j Z(x_j)+\gamma,$$
for $n$ large enough, where $|A|$ denotes the Lebesgue measure of $A$. Therefore letting $n\to\infty$ we have
$$ \inf_{B_{\rho_1(x_i)}} W(x)\leq Z(x_i)-\min_j Z(x_j)+\gamma,$$
and hence
$$ W(x_i)\leq Z(x_i)-\min_j Z(x_j)+2\gamma.$$
Similarly using the upper bound in Theorem \ref{invmeas} we have
$$ W(x_i)\geq Z(x_i)-\min_j Z(x_j)-2\gamma.$$
Since $\gamma>0$ is arbitrary, we have $W(x_i)= Z(x_i)-\min_j Z(x_j)$. This completes the proof of the theorem.\hfill $\Box$.

\spa

\noindent {\bf Other minor corrections}\\
\begin{itemize}
 \item The closed bounded set $\Gamma$ in the expression (15), page 479 in \cite{bisbor}, should be taken as a closed bounded ball with center at the origin.

\item We need to assume that $\frac{\partial a_{ij}}{\partial x}(0)=0$, for all $i, j$, in Theorem 3 in \cite{bisbor}.

\end{itemize}

\noindent \textbf{\large Appendix}\\

Here we prove uniqueness of viscosity solution of the Hamilton Jacobi equation
\begin{eqnarray}
 u_t + H(x, Du) &=& 0\ \ \mbox{in}\ \ \R^d\times (0, T]\label{unihjb}
\\
u(x, 0)&=& g(x).\nonumber
\end{eqnarray}
We prove the uniqueness along the lines of \cite{crandall}. We assume the following
\begin{itemize}
 \item[(A1)] For some local modulus $\sigma : [0,\infty)\times [0, \infty)\to [0, \infty)$
$$H(y, \lambda(x-y))-H(x,\lambda(x-y))\leq \sigma (\lambda \|x-y\|^2+\|x-y\|, R)$$
for $R>0, x, y\in B_R, \lambda \geq 0$.
\item [(A2)] There exists a constant $C_0>0$ such that
$$ |H(x, p)-H(x, q)|\leq C_0(1+\|p\|+\|q\|)\|p-q\|$$
for all $x, p, q\in\R^d$.
\end{itemize}
Let $u , v$ be two viscocity solution of (\ref{unihjb}) which are continuous on $\R^d\times[0, T]$ and Lipschitz in space variable uniformly w.r.t.\ $t\in[0, T]$.
To show that $u=v$ on $\R^d\times[0, T]$ it is enough to show that
$$ \sup_{\R^d\times[0, T]}(u-v)^+\leq 0.$$

\noi Now if possible, let there exists $\chi>0$ such that for some $t_0>0$, $x_0\in\R^d$,
$$ u(x_0,t_0)-v(x_0,t_0)\geq 3\chi.$$
Let
$$ |u(x,t)-u(y,t)|\leq K \|x-y\|\ \mbox{and}\ |v(x,t)-v(y,t)|\leq K \|x-y\|$$
for all $t\in[0, T]$. Therefore from assumption (A2) we have
\begin{eqnarray}
|H(x,p)-H(x,q)|\leq C\|p-q\|\label{bb123}
\end{eqnarray}
for $p\in D^+_x u(x,t),\ q\in D^-_y v(y,s)$ for all $x,y\in\R^d$ and $s,t\in(0,T]$. For $0<\theta<1$, define
$$ w_\delta(x,t)=\exp{\left\{\frac{1}{\delta}\Big(\log(1+\|x\|^2)+\frac{CT^\theta}{1-\theta}t^{1-\theta}-R\Big)\right\}}.$$
It is easy to check that
$$ \frac{\partial w_\delta}{\partial t} - C\frac{T^\theta}{t^\theta}\|Dw_\delta\|\geq 0 \ \mbox{in}\ \R^d\times(0, \infty).$$
Choose $R$ large enough so that $\log(1+\|x_0\|^2)+\frac{CT^\theta}{1-\theta}t_0^{1-\theta}-R<0$. Now define $Q_R=\{\delta\log w_\delta(x,t)\leq 1\}\subset \R^d\times[0, T]$.
It is easy to see that there exists a compact set $\mathbb{K}$ independent of $\delta$ such that $Q_R\subset\mathbb{K}\times[0, T]$ for all $\delta>0$. Let
$M=(\max_{Q_R}u(x,t)\vee\max_{Q_R}v(x,t))$. Now we choose $\delta$ small so that $w_\delta(x_0,t_0)<(\chi\wedge 1)$ and $e^{\frac{1}{2\delta}}>2M$. On $Q_R\times Q_R$,
define
$$\Phi(x,t,y,s)=u(x,t)-v(y,s)-w_\delta(x,t)-\frac{\|x-y\|^2+|t-s|^2}{\eps}-a(t+s)$$
for $a, \eps>0$. For $a>0$ small, we have
$$ \sup_{Q_R\times Q_R}\Phi \geq \Phi(x_0, t_0, x_0, t_0)\geq \chi>0.$$
Let $(\bar x,\bar t, \bar y, \bar s)$ be the point where the maximum of $\Phi$ is attended. By standard calculations, it can be shown that
\begin{itemize}
 \item $\frac{\|\bar x-\bar y\|^2+|\bar t-\bar s|^2}{\eps}\leq 2M$.
\item $\frac{\|\bar x-\bar y\|^2}{\eps}\to 0$ as $\eps\to 0$.
\end{itemize}
Using the uniform continuity of $u, v, w_\delta$ on $Q_R$ and by our choice of $\delta$ it is easy to show that $\bar t>0, \bar s>0$ for $\eps>0$ small.
Also for the same reasons, $\delta\log w_\delta(\bar x, \bar t)<1$ and $\delta\log w_\delta(\bar y, \bar s)<1$. Hence we can apply the definition of viscosity
solution to the points $(\bar x, \bar t)$, $(\bar y, \bar s)$ to get
$$ a + \frac{\partial w_\delta(\bar x,\bar t)}{\partial t} + 2\frac{\bar t-\bar s}{\eps} + H\left(\bar x, Dw_\delta(\bar x,\bar t)+2\frac{\bar x-\bar y}{\eps}\right)\leq 0,$$
and
$$ -a + 2\frac{\bar t-\bar s}{\eps} + H\left(\bar y, 2\frac{\bar x-\bar y}{\eps}\right)\geq 0.$$

Subtracting we have
\begin{eqnarray*}
2a + \frac{\partial w_\delta(\bar x,\bar t)}{\partial t} + H\left(\bar x, Dw_\delta(\bar x,\bar t)+2\frac{\bar x-\bar y}{\eps}\right)
- H\left(\bar x,2\frac{\bar x-\bar y}{\eps}\right)&&
\\
 + H\left(\bar x,2\frac{\bar x-\bar y}{\eps}\right) - H\left(\bar y,2\frac{\bar x-\bar y}{\eps}\right) \leq  0 &&
\end{eqnarray*}
Therefore using (\ref{bb123}) we have
\begin{eqnarray*}
2a + \frac{\partial w_\delta(\bar x,\bar t)}{\partial t} - \frac{CT^\theta}{\bar t^\theta}|Dw_\delta(\bar x,\bar t)|
 + H\left(\bar x,2\frac{\bar x-\bar y}{\eps}\right) - H\left(\bar y,2\frac{\bar x-\bar y}{\eps}\right) \leq  0
\end{eqnarray*}
Now letting $\eps\to 0$ and using (A1) we have $2a\leq 0$ which is a contradiction. This proves the uniqueness claim.